\documentclass[12pt,leqno]{article}

\usepackage{amsmath,amssymb,amsthm}

\makeatletter

\newtheorem{Theorem}{Theorem}

\newtheorem{Lemma}[Theorem]{Lemma}
\newtheorem{Corollary}[Theorem]{Corollary}

\theoremstyle{definition}

\theoremstyle{remark}

\def\({{\rm (}}
\def\){{\rm )}}

\let\Mathrm\operator@font
\let\Cal\mathcal

\def\standop#1{\mathop{\Mathrm #1}\nolimits}
\def\difstop#1#2{\expandafter\def\csname #1\endcsname{\standop{#2}}}
\def\defstop#1{\difstop{#1}{#1}}

\defstop{Hom}
\defstop{Gal}
\defstop{Spec}

\defstop{Ker}
\defstop{Min}
\difstop{Image}{Im}
\difstop{tdeg}{trans.deg}

\defstop{length}
\defstop{supp}
\defstop{Proj}
\defstop{Flat}
\difstop{height}{ht}

\def\specialarrow#1{\setbox\z@=\hbox{$\m@th
 \mathop{\vphantom{\rightarrow}}\limits^{\hspace{.5ex}{#1}\hspace
{.8ex}}$}\mathrel{\ifdim\wd\z@<1.2em\dimen\tw@
1.2em\else\dimen\tw@\wd\z@\fi\copy\z@\kern-\wd\z@\hbox to\dimen\tw@
{\rightarrowfill}}}

\def\sdarrow#1{\downarrow\hbox to 0pt{\scriptsize$#1$\hss}}
\def\suarrow#1{\uparrow\hbox to 0pt{\scriptsize$#1$\hss}}

\def\section{\@startsection{section}{1}{\z@ }%
{-3.5ex plus -1ex minus -.2ex}{2.3ex plus .2ex}{\bf }}

\long\def\refname{\par\kern -3ex
\begin{center}\rm R\sc{eferences}\end{center}\par\kern 
-2ex}

\def\@seccntformat#1{\csname the#1\endcsname.\quad}

\def\@@@sect#1#2#3#4#5#6[#7]#8{%
   \ifnum #2>\c@secnumdepth 
      \def \@svsec {}\else \refstepcounter {#1}%
      \def\@svsec{}
   \fi 
   \@tempskipa #5\relax 
   \ifdim \@tempskipa >\z@ 
     \begingroup #6\relax \@hangfrom {\hskip #3\relax 
     \@svsec}{\interlinepenalty \@M #8\par }\endgroup 
     \csname #1mark\endcsname {#7}
   \else 
   \def \@svsechd {#6\hskip #3\@svsec #8\csname #1mark\endcsname {#7}}
   \fi \@xsect {#5}}

\def\@@@startsection#1#2#3#4#5#6{%
 \if@noskipsec \leavevmode \fi \par \@tempskipa #4\relax \@afterindenttrue 
 \ifdim \@tempskipa <\z@ \@tempskipa -\@tempskipa \@afterindentfalse 
 \fi \if@nobreak \everypar {}\else \addpenalty {\@secpenalty }\addvspace 
  {\@tempskipa }\fi \@ifstar {\@ssect {#3}{#4}{#5}{#6}}{\@dblarg 
  {\@@@sect {#1}{#2}{#3}{#4}{#5}{#6}}}}

\def\theparagraph{\thesection.\arabic{paragraph}}
\def\aparagraph{\@@@startsection{paragraph}{2}{\z@ }%
              {1.75ex plus .2ex minus .15ex}{-1em}{\bf(\theparagraph) } }
\def\paragraph{\@@@startsection{paragraph}{2}{\z@ }%
              {1.75ex plus .2ex minus .15ex}{-1em}{}{\bf(\theparagraph)} }

\let\c@Theorem\c@paragraph

\title{Another proof of global $F$-regularity of Schubert varieties\thanks
{2000 Mathematics Subject Classification. Primary 14M15, Secondary 13A35.}}

\author{M{\sc itsuyasu} H{\sc ashimoto}}
\date{}

\begin{document}

\maketitle

\begin{abstract}
Recently, Lauritzen, Raben-Pedersen and Thomsen proved 
that Schubert varieties are globally $F$-regular.
We give another proof.
\end{abstract}

\section{Introduction}

Let $p$ be a prime number, $k$ an algebraically closed field of 
characteristic $p$, and $G$ a simply connected semisimple affine 
algebraic group over $k$.
Let $T$ be a maximal torus of $G$.
We choose a base of the root system of $G$.
Let $B$ be the negative Borel subgroup of $G$.
Let $P$ be a parabolic subgroup of $G$ containing $B$.
The closure of a $B$-orbit on $G/P$ is called a Schubert variety.

Recently, Lauritzen, Raben-Pedersen and Thomsen proved 
that Schubert varieties are globally $F$-regular \cite{LRT} 
utilizing Bott-Samelson resolution.
The objective of this paper is to give another proof.
Our proof depends on a simple inductive argument utilizing the familiar
technique of fibering the Schubert variety as a $\Bbb P^1$-bundle
over a smaller Schubert variety.

Global $F$-regularity was first defined by Smith \cite{Smith}.
A projective variety over $k$ is said to be globally $F$-regular
if it admits a strongly $F$-regular homogeneous coordinate ring.
As a corollary, we have that the all local rings of a Schubert variety
is $F$-regular, in particular, $F$-rational, Cohen-Macaulay and normal.

A globally $F$-regular variety is Frobenius split.
It has long been known that Schubert varieties are Frobenius split
\cite{MR}.
Given an ample line bundle over $G/P$, the associated projective
embedding of a Schubert variety of $G/P$ is projectively normal 
\cite{RR} and arithmetically Cohen-Macaulay \cite{Ramanathan}.
We can prove that the coordinate ring is strongly $F$-regular in fact.

Over globally $F$-regular varieties, there are some nice vanishing theorems.
One of these gives a short proof of Demazure's vanishing theorem.

Acknowledgement. The author is grateful to Professor V. B. Mehta for
valuable advice.
In particular, Corollary~\ref{Mehta.thm} is due to him.
He also kindly 
showed the result of  Lauritzen, Raben-Pedersen and Thomsen to the author.
Special thanks are also due to Professor V. Srinivas and K.-i.~Watanabe 
for valuable advice.

\section{Preliminaries}

Let $p$ be a prime number, and $k$ an algebraically closed field of 
characteristic $p$.
For a ring $A$ of characteristic $p$, the Frobenius map $A\rightarrow A$ 
$(a\mapsto a^p)$ is denoted by $F$ or $F_A$.
So $F_A^e$ maps $a$ to $a^{p^e}$ for $a\in A$ and $e\geq 0$.

Let $A$ be a $k$-algebra.
The ring $A$ with the $k$-algebra structure given by
\[
k\xrightarrow{F^{-r}_k} k \rightarrow A
\]
is denoted by $A^{(r)}$ for $r\in\Bbb Z$.
Note that $F_A^e\colon A^{(r+e)}\rightarrow A^{(r)}$ is a $k$-algebra map
for $e\geq 0$ and $r\in\Bbb Z$.
For $a\in A$ and $r\in \Bbb Z$, the element $a$ viewed as an element in 
$A^{(r)}$ is sometimes denoted by $a^{(r)}$.
So $F_A^e(a^{(r+e)})=(a^{(r)})^{p^e}$ for $a\in A$, $r\in\Bbb Z$ and $e\geq 0$.

Similarly, for a $k$-scheme $X$ and $r\in\Bbb Z$, the $k$-scheme $X^{(r)}$ 
is defined.
The Frobenius morphism $F_X^e\colon X^{(r)}\rightarrow X^{(r+e)}$ is a 
$k$-morphism.

A $k$-algebra $A$ is said to be $F$-finite if the Frobenius map
$F_A\colon A^{(1)}\rightarrow A$ is finite.
A $k$-scheme $X$ is said to be $F$-finite if the Frobenius morphism
$F_X\colon X\rightarrow X^{(1)}$ is finite.
Let $A$ be an $F$-finite Noetherian $k$-algebra.
We say that $A$ is strongly $F$-regular if for
any non-zerodivisor $c\in A$, there exists some $e\geq 0$ such that
$cF^e_A\colon A^{(e)}\rightarrow A$ $(a^{(e)}\mapsto ca^{p^e})$ is a
split monomorphism as an $A^{(e)}$-linear map \cite{HH}.
A strongly $F$-regular $F$-finite ring is $F$-rational in the sense of
Fedder--Watanabe \cite{FW}, and is Cohen--Macaulay normal.

Let $X$ be a quasi-projective $k$-variety.
We say that $X$ is globally $F$-regular if 
for any invertible sheaf $\Cal L$ over $X$ and 
any $a\in\Gamma(X,\Cal L)\setminus 0$, 
the composite
\[
\Cal O_{X^{(e)}}\rightarrow F_*^e \Cal O_X \xrightarrow{F^e_* a} 
F_*^e\Cal L
\]
has an $\Cal O_{X^{(e)}}$-linear splitting \cite{Smith}, \cite{Hashimoto}.
$X$ is said to be $F$-regular if $\Cal O_{X,x}$
is strongly $F$-regular for any closed point $x$ of $X$.

Smith \cite[(3.10)]{Smith} proved the following fundamental theorem on
global $F$-regularity.
See also \cite[(3.4)]{Watanabe} and \cite[(2.6)]{Hashimoto}.

\begin{Theorem}\label{Smith.thm}
Let $X$ be a projective variety over $k$.
Then the following are equivalent.
\begin{description}
\item[1] There exists some ample Cartier divisor $D$ on $X$ such that the
section ring $\bigoplus_{n\in\Bbb Z}\Gamma(X,\Cal O(nD))$ is strongly
$F$-regular.
\item[2] The section ring of $X$ with respect to every ample 
Cartier divisor is strongly $F$-regular.
\item[3] There exists some ample effective Cartier divisor $D$ on $X$ such 
that there exists some $e\geq 0$ and an $\Cal O_{X^{(e)}}$-linear splitting
of $\Cal O_{X^{(e)}}\rightarrow F_*^e \Cal O_X \rightarrow F_*^e \Cal O(D)$
and that the open set $X-D$ is $F$-regular.
\item[4] $X$ is globally $F$-regular.
\end{description}
\end{Theorem}

A globally $F$-regular variety is $F$-regular.

An affine $k$-variety $\Spec A$ is globally $F$-regular
if and only if $A$ is strongly $F$-regular if and only if $\Spec A$ is
$F$-regular.

A globally $F$-regular variety is Frobenius split in the sense of 
Mehta-Ramanathan \cite{MR}.
As the theorem above shows, if $X$ is a globally $F$-regular projective 
variety, then the section
ring of $X$ with respect to every ample divisor is Cohen--Macaulay normal.


The following is a useful lemma.

\def\citinfo{\cite[Proposition~1.2]{HWY}}
\begin{Lemma}[\citinfo]\label{easy.thm}
Let $f\colon X\rightarrow Y$ be a $k$-morphism between projective 
$k$-varieties.
If $X$ is globally $F$-regular and $\Cal O_Y\rightarrow f_*\Cal O_X$ is
an isomorphism, then $Y$ is globally $F$-regular.
\end{Lemma}

Let $G$ be a simply connected semisimple algebraic group over $k$,
and $T$ a maximal torus of $G$.
We fix a base of the set of roots of $G$.
Let $B$ be the negative Borel subgroup.
Let $P$ be a parabolic subgroup of $G$ containing $B$.
Then $B$ acts on $G/P$ from the left.
The closure of a $B$-orbit of $G/P$ is called a Schubert variety.
Any $B$-invariant closed subvariety of $G/P$ is a Schubert variety.
The set of Schubert varieties in $G/B$ and the Weyl group $W(G)$ of $G$ 
are in one-to-one correspondence.
For a Schubert variety $X$ in $G/B$, there is a unique $w\in W(G)$ such that 
$X=\overline{BwB/B}$, where the overline denotes the closure operation.
We need the following theorem
later.

\begin{Theorem}\label{normal.thm}
A Schubert variety in $G/P$ is a normal variety.
\end{Theorem}

For the proof, see \cite[Theorem~3]{RR}, \cite{Andersen}, \cite{Seshadri},
and \cite{MS}.

Let $X$ be a Schubert variety in $G/P$.
Then $\tilde X=\pi^{-1}(X)$ is a $B$-invariant reduced subscheme of $G/B$,
where $\pi\colon G/B\rightarrow G/P$ is the canonical projection.
It has a dense $B$-orbit, and actually $\tilde X$ 
is a Schubert variety in $G/B$.

Let $Y=\rho^{-1}(\tilde X)$, where $\rho\colon G\rightarrow G/P$ is the
canonical projection.
Let $\Phi\colon Y\times P/B \rightarrow Y\times_{X}\tilde X$ be the 
$Y$-morphism given by $\Phi(y,pB)=(y,ypB)$.
Since $(y,\tilde x B)\mapsto (y,y^{-1}\tilde x B)$ gives the inverse, 
$\Phi$ is an isomorphism.
Note that $(p_1)_* \Cal O_{Y\times P/B}\cong \Cal O_Y$, where $p_1\colon
Y\times P/B\rightarrow Y$ is the first projection, since 
$P/B$ is a $k$-complete variety and $H^0(P/B,\Cal O_{P/B})=k$.
As $\Phi$ is a $Y$-isomorphism, we have that 
$(\pi_1)_* \Cal O_{Y\times_X \tilde X}\cong \Cal O_Y$,
where $\pi_1\colon Y\times_X \tilde X\rightarrow Y$ is the first
projection.
As $\pi_1$ is a base change of $\pi\colon \tilde X\rightarrow X$ by 
the faithfully flat morphism $Y\rightarrow X$, we have

\begin{Lemma}\label{p-to-b.thm}
$\pi_*\Cal O_{\tilde X}\cong \Cal O_X$.
In particular, if $\tilde X$ is globally $F$-regular, then so is $X$.
\end{Lemma}

Let $w\in W(G)$, and $X=X_w$ be the corresponding Schubert variety 
$\overline{BwB/B}$ in $G/B$.
Assume that $w$ is nontrivial.
Then there exists some simple root $\alpha$ such that $l(w s_\alpha)
=l(w)-1$, where $s_\alpha$ is the reflection corresponding to $\alpha$, and
$l$ denotes the length.
Set $X'=X_{w'}$ be the Schubert variety $\overline{Bw'B/B}$, where 
$w'=ws_\alpha$.
Let $P_\alpha$ be the minimal parabolic subgroup $Bs_\alpha B\cup B$.
Let $Y$ be the Schubert variety $\overline{BwP_\alpha/P_\alpha}$.

The following is due to Kempf \cite[Lemma~1]{Kempf}.

\begin{Lemma}\label{Kempf.thm}
Let $\pi_\alpha\colon G/B\rightarrow G/P_\alpha$ be the
canonical projection.
Then $X'$ is birationally mapped onto $Y$.
In particular, $(\pi_\alpha)_* \Cal O_{X'}=\Cal O_Y$ 
(by Theorem~\ref{normal.thm}).
We have $(\pi_\alpha)^{-1}(Y)=X$, and $\pi|_X\colon X\rightarrow Y$ is 
a $\Bbb P^1$-fibration, hence is smooth.
\end{Lemma}

Let $X$ be a Schubert variety in $G/B$.
Let $\rho$ be the half-sum of positive roots, and set $\Cal L=
\Cal L((p-1)\rho)|_X$, where $\Cal L((p-1)\rho)$ is the invertible sheaf
on $G/B$ corresponding to the weight $(p-1)\rho$.
The following was proved by Ramanan--Ramanathan \cite{RR}.
See also Kaneda \cite{Kaneda}.

\begin{Theorem}\label{rr.thm}
There is a section $s\in H^0(X,\Cal L)\setminus 0$ such that
the composite 
\[
\Cal O_{X^{(1)}}\rightarrow F_* \Cal O_X \xrightarrow {F_* s} F_*\Cal L
\]
splits.
\end{Theorem}

Since $\Cal L$ is ample, we immediately have the following.

\begin{Corollary}\label{Mehta.thm}
$X$ is globally $F$-regular if and only if $X$ is $F$-regular.
\end{Corollary}

\proof The `only if' part is obvious.
The `if' part follows from the theorem and Theorem~\ref{Smith.thm}, {\bf 
3$\Rightarrow$4}.
\qed

\section{Main theorem}

Let $k$ be an algebraically closed field, $G$ a semisimple simply connected
algebraic group over $k$, $T$ a maximal torus of $G$.
We fix a basis of the set of roots of $G$, and let $B$ be
the negative Borel subgroup of $G$.

In this section we prove the following theorem.

\begin{Theorem}\label{main.thm}
Let $P$ be a parabolic subgroup of $G$ containing $B$, and
let $X$ be a Schubert variety in $G/P$.
Then $X$ is globally $F$-regular.
\end{Theorem}

\proof Let $\pi\colon G/B\rightarrow G/P$ be the canonical projection, and set 
$\tilde X=\pi^{-1}(X)$.
Then $\tilde X$ is a Schubert variety in $G/B$.
By Lemma~\ref{p-to-b.thm}, it suffices to show that $\tilde X$ is 
globally $F$-regular.
So in the proof, we may and shall assume that $P=B$.

So let $X=\overline{BwB/B}$.
We proceed by induction on the dimension of $X$, in other words, $l(w)$.
If $l(w)=0$, then $X$ is a point and $X$ is globally $F$-regular.
Let $l(w)>0$.
Then there exists some simple root $\alpha$ such that $l(ws_\alpha)=l(w)-1$.
Set $w'=ws_\alpha$, $X'=\overline{Bw'B/B}$, $P_\alpha=Bs_\alpha B\cup B$, 
and $Y=\overline{Bw P_\alpha/P_\alpha}$.

By induction assumption, $X'$ is globally $F$-regular.
By Lemma~\ref{Kempf.thm} and Lemma~\ref{easy.thm}, $Y$ is also globally
$F$-regular.
In particular, $Y$ is $F$-regular.
By Lemma~\ref{Kempf.thm}, $X\rightarrow Y$ is smooth.
By \cite[(4.1)]{LS}, $X$ is $F$-regular.
By Corollary~\ref{Mehta.thm}, $X$ is globally $F$-regular.
\qed

\def\citinfo{Demazure's vanishing \cite{RR}, \cite{Kaneda}}
\begin{Corollary}[\citinfo]
Let $X$ be a Schubert variety in $G/B$, $\lambda$ a dominant weight, and
$\Cal L:=\Cal L(\lambda)|_X$.
Then $H^i(X,\Cal L)=0$ for $i>0$.
\end{Corollary}

\proof This follows from the theorem and \cite[(4.3)]{Smith}.
\qed

Let $P$ be a parabolic subgroup of $G$ containing $B$.
Let $X$ be a Schubert variety in $G/P$.
Let $\Cal M_1,\ldots,\Cal M_r$ be effective line bundles on $G/P$,
and set $\Cal L_i:=\Cal M_i|_X$.
In \cite{KR}, Kempf and Ramanathan proved that the $k$-algebra
$C:=\bigoplus_{\mu\in\Bbb N^r} \Gamma(X,\Cal L_\mu)$ has rational 
singularities, where $\Cal L_\mu=\Cal L_1^{\otimes \mu_1}\otimes\cdots
\otimes \Cal L_r^{\otimes \mu_r}$ for $\mu=(\mu_1,\ldots,\mu_r)\in\Bbb Z^r$.
We can prove a very similar result.

\begin{Corollary}
Let the notation be as above.
The $k$-algebra $C$ is strongly $F$-regular.
\end{Corollary}

By \cite[Theorem~2.6]{Hashimoto}, $\tilde C=\bigoplus_{\mu\in\Bbb Z^r}
\Gamma(X,\Cal L_\mu)$ is a quasi-$F$-regular domain.
By \cite[Lemma~2.4]{Hashimoto}, $C$ is also quasi-$F$-regular.
By \cite[Theorem~2]{RR}, $C$ is finitely generated over $k$, and
$C$ is strongly $F$-regular.
\qed


\begin{flushleft}
Graduate School of Mathematics\\
Nagoya University\\
Chikusa-ku,  Nagoya 464--8602\\
JAPAN
\end{flushleft}

\begin{flushleft}
\small {\em E-mail address}: \tt hasimoto@math.nagoya-u.ac.jp
\end{flushleft}

\end{document}